\documentclass[12pt,twoside,reqno]{amsart}
\usepackage[colorlinks=true,citecolor=blue]{hyperref}
\usepackage{mathptmx, amsmath, amssymb, amsfonts, amsthm, mathptmx, enumerate, color,mathrsfs}
\usepackage{graphicx}
\usepackage{float}  

\usepackage{epstopdf}

\usepackage{booktabs}
\usepackage{tabularx}
\usepackage{caption}
\usepackage{diagbox}
\usepackage{multirow}
\usepackage{longtable}
\usepackage{threeparttable}
\usepackage{booktabs} 
\usepackage{array} 
\usepackage{paralist} 
\usepackage{verbatim} 
\usepackage{subfig} 
\usepackage{amsmath}
\usepackage{amssymb}

\theoremstyle{definition}

\numberwithin{equation}{section}

\begin{document}
\setcounter{page}{1}

\vspace*{2.0cm}
\title[Research on the Evaluation Index System]
{Research on the Evaluation Index System of Enterprise Production Efficiency}
\author[ W. Li,J. Cai,C. Wang, Y.Chen,J.Xu, J.Zhao,Y.Chen]{  WEN LI$^1$,JIAXI CAI$^1$,CAIPING WANG$^2$, YANGCE CHEN$^1$,JIN XU$^3$, JUNJIAN ZHAO$^1$,YASONG CHEN $^{1,*}$}
\maketitle
\vspace*{-0.6cm}

\begin{center}
{\footnotesize

$^1$School of Mathematical Sciences, Tiangong University, Tianjin300387, China\\
$^2$Zhangjiakou Cigarette Factory Co., Ltd,Zhangjiakou075000, China\\
$^3$Ximatou Central School, Daliuhe Town, Wen'an County, Langfang065800,China\\
}\end{center}

\vskip 4mm {\footnotesize \noindent {\bf Abstract.}
This paper focuses on studying the evaluation index system for the production efficiency of tobacco enterprises. Considering the limitations of existing evaluation methods in accurately assessing the production quality of cigarette enterprises, a mathematical model based on the Analytic Hierarchy Process (AHP) is established. This model constructs an evaluation framework for the production efficiency of cigarette enterprises and subsequently analyzes the significance of each index within this framework. To comprehensively analyze the multi-index and feasibility aspects of the selected projects, the AHP method is employed to establish a comprehensive feasibility research and evaluation structure model. The result of this feasibility study provides the conclusion that the construction of an evaluation index system for the production efficiency of cigarette enterprises can indeed promote the enhancement of their production efficiency.



\renewcommand{\thefootnote}{}
\footnotetext{ $^*$Corresponding author.
\par
E-mail addresses:chenyasong@tiangong.edu.cn (Y. Chen).
\par
Received xxxx xx, 2024; Accepted xxxx xx, 2024. }

\section{Introduction}
Taking comprehensive production efficiency as a comprehensive performance indicator, applying it to track and evaluate the production capacity of enterprises, and providing critical support for monitoring, evaluating, alerting, and analyzing influencing factors of production efficiency. With the continuous growth and development of China's cigarette industry, conducting more efficient product manufacturing and quality management becomes even more crucial\cite{0.5}. For the cigarette industry, although some enterprises have explored the construction of production efficiency systems, the practical results have not effectively promoted the improvement of production efficiency\cite{0.6}. Utilizing comprehensive production efficiency as a comprehensive performance indicator, applied to track and evaluate the production capacity and efficiency of manufacturing enterprises, plays a crucial supportive role in monitoring, evaluating, alerting, and analyzing factors influencing production efficiency\cite{1}.

As a typical process-oriented manufacturing enterprise, cigarette companies operate in a complex and diverse production environment, influenced by various factors affecting production efficiency. These factors include the ratio of good products to total output (quality factor), the ratio of production speed to design speed (performance factor), the ratio of actual normal operating time to planned production time (availability factor), the ratio of production alignment to planned production time (production synergy factor), and the ratio of equipment failure downtime to planned production time (equipment failure factor), among others\cite{1.4}.

In the past, there have been limitations in the collection, storage, and analysis of data. When analyzing and processing relevant information, companies often have to extract and analyze the limited information they can obtain to the maximum extent, which inadvertently increases their workload. Moreover, issues such as incomplete and outdated information directly affect the analysis results. Therefore, manual statistical analysis is time-consuming, labor-intensive, and unable to achieve efficient and scientific analysis\cite{1.5}. Although some companies have achieved the calculation and visualization of production efficiency indicators through information means, they are limited by data acquisition capabilities, analysis tools, and computing power, resulting in coarse-grained analysis, insufficient analysis basis, significant analysis lag, and distorted analysis results, which can easily mislead decision-making\cite{2}.

Therefore, it is necessary to construct a system and conduct research on production efficiency indicators for cigarette enterprises, for the reference of the majority of cigarette industry enterprises.

\section{Preliminaries}

\subsection{The principles for establishing evaluation indicators}

For complex systems like cigarette manufacturing efficiency in enterprises, it is necessary to use multiple indicators to form an organic whole to describe its state and changes. It is currently not possible to accurately assess it using only a few indicators\cite{3}.

As a whole composed of multiple indicators, there exists a certain connection and interaction among the various indicators in the evaluation indicator system, which links the evaluators and the objects being evaluated. A reasonable evaluation indicator system is the cornerstone for comprehensive evaluation of any object, and without it, accurate assessment cannot be achieved. Therefore, the scientificity of the indicator system directly affects the rationality of the evaluation results\cite{3.5}.

The selection of evaluation indicator system must adhere to the following principles\cite{4}:

(1)Principle of Scientific: The selection of evaluation indicators should first and foremost possess strong scientific validity. The chosen indicators must be in line with the research content of the article and meet practical requirements.

(2)Principle of Operability: The selection of evaluation indicators should be practical and feasible, combining theoretical analysis with practical application. It is important to choose indicators that have a significant impact and for which relevant data can be easily obtained.

(3)Principle of Holism: Evaluation is a comprehensive and systematic process, so when selecting indicators, it is important to consider their systemic and holistic nature. It is necessary to continuously refine and revise the indicators, while selecting as many indicators as possible that cover various subsystems of the evaluation object to meet the requirements of system development.

(4)Principle of Optimization: The selection of evaluation indicators should maximize the inclusion of various aspects of the research content, but it should not be overly complex or redundant. Therefore, when selecting evaluation indicators, the focus should be on selecting representative key indicators from various aspects, aiming for a small yet precise set of indicators, in accordance with the principle of optimization.

\subsection{The constructed indicator system in this article}

After reviewing a large number of relevant literature, the authors of this article found that the existing indicator systems for evaluating cigarette companies mostly focus on individual technical evaluations and do not establish a comprehensive efficiency evaluation system. Based on the integration of previous research findings and consultation with relevant experts, this article initially determined the secondary and tertiary indicators\cite{4.5}. Through brainstorming with a certain cigarette enterprises experts, the bloated indicators were reduced and the neglected indicators were added\cite{5}.

The comprehensive evaluation index system for cigarette enterprises efficiency is constructed based on the production and construction process of cigarette enterprises. It consists of three levels: goal level $A$, criterion level $B$, and index level $C$.

The goal level A is the efficiency evaluation index system for cigarette enterprises.

 The criterion level B includes 6 indicators: plan execution, equipment efficiency, production material consumption, process control of silk production quality, process control of rolling and packaging quality, energy saving and emission reduction.

 The criterion level C includes 23 indicators: production plan completion rate, unit output maintenance duration, unit output maintenance frequency, rolling equipment operating efficiency, packaging equipment operating efficiency, leaf-to-silk ratio, consumption of cigarette paper per carton, consumption of small box paper per carton, residual cigarettes per carton, consumption of cigarette filters per carton, consumption of cut cigarette per carton, number of process deviations, acceptance rate of premium quality finished products, moisture deviation in thin sheet drying machine export, moisture deviation in fuel drying machine export, $CPK$qualification rate of cigarette quality, $CPK$qualification rate of draw resistance quality, deviation in single cigarette weight, score of finished product release inspection, premium product rate, energy consumption per carton, pollution emissions per carton, carbon dioxide emissions per carton.

As shown in the Table \ref{Table zhibiaotixi}:

\begin{table}[htbp]
\centering
\caption{Evaluation Index System for Manufacturing Efficiency in Cigarette Enterprises}
\label{Table zhibiaotixi}
\resizebox{\linewidth}{!}{
\begin{tabular}{ccc}
\toprule
 Goal Level $A$      &  Criterion Level $B$    & Index Level $C$  \\
\midrule 
\multirow{23}{*}{$A$  $\mathpunct{:}$ Efficiency Evaluation Index System for Cigarette Enterprises} &   \multirow{1}{*}{$B_{1}$ $\mathpunct{:}$ Plan Execution} & $C_{11}$ $\mathpunct{:}$ Production Plan Completion Rate \\
\cline{2-3}
 &  \multirow{4}{*}{$B_{2}$ $\mathpunct{:}$ Equipment Efficiency} & $C_{21}$ $\mathpunct{:}$ Unit Output Maintenance Duration  \\
 &   & $C_{22}$ $\mathpunct{:}$ Unit Output Maintenance Frequency  \\
 &   & $C_{23}$ $\mathpunct{:}$ Rolling Equipment Operating Efficiency  \\
 &   & $C_{24}$ $\mathpunct{:}$ Packaging Equipment Operating Efficiency  \\
\cline{2-3}
 &  \multirow{6}{*}{$B_{3}$ $\mathpunct{:}$ Production Material Consumption}  & $C_{31}$ $\mathpunct{:}$ Leaf-to-Silk Ratio  \\
 &   &  $C_{32}$ $\mathpunct{:}$ Consumption of Cigarette Paper Per Carton  \\
 &   &  $C_{33}$ $\mathpunct{:}$ Consumption of Small Box Paper Per Carton  \\
 &   &  $C_{34}$ $\mathpunct{:}$ Residual Cigarettes Per Carton  \\
 &   &  $C_{35}$ $\mathpunct{:}$ Consumption of Cigarette Filters Per Carton  \\
 &   &  $C_{36}$ $\mathpunct{:}$ Consumption of Cut Cigarette Per Carton  \\
\cline{2-3}
 &  \multirow{4}{*}{$B_{4}$ $\mathpunct{:}$ Process Control of Silk Production Quality}  & $C_{41}$ $\mathpunct{:}$ Number of Process Deviations  \\
 &   & $C_{42}$ $\mathpunct{:}$ Acceptance Rate of Premium Quality Finished Products  \\
 &   & $C_{43}$ $\mathpunct{:}$ Moisture Deviation in Thin Sheet Drying Machine Export \\
 &   & $C_{44}$ $\mathpunct{:}$ Moisture Deviation in Fuel Drying Machine Export \\
\cline{2-3}
 &  \multirow{5}{*}{$B_{5}$ $\mathpunct{:}$ Process Control of Rolling and Packaging Quality}  & $C_{51}$ $\mathpunct{:}$ $CPK$Qualification Rate of Cigarette Quality  \\
 &   & $C_{52}$ $\mathpunct{:}$ $CPK$Qualification Rate of Draw Resistance Quality    \\
 &   & $C_{53}$ $\mathpunct{:}$ Deviation in Single Cigarette Weight  \\
 &   & $C_{54}$ $\mathpunct{:}$ Score of Finished Product Release Inspection  \\
 &   & $C_{55}$ $\mathpunct{:}$ Premium Product Rate  \\
\cline{2-3}
 &  \multirow{3}{*}{$B_{6}$ $\mathpunct{:}$ Energy Saving and Emission Reduction}  & $C_{61}$ $\mathpunct{:}$ Energy Consumption Per Carton  \\
 &   &  $C_{62}$ $\mathpunct{:}$ Pollution Emissions Per Carton  \\
 &   &  $C_{63}$ $\mathpunct{:}$ Carbon Dioxide Emissions Per Carton  \\
\bottomrule 
\end{tabular}
}
\end{table}

\section{Basic concept}

The Analytic Hierarchy Process (AHP) was proposed by American operations researcher Professor Thomas L. Saaty at the University of Pittsburgh in the 1970s\cite{6}. This method primarily divides the influencing factors of a research problem into multiple levels and combines quantitative and qualitative analysis methods to further analyze them. The basic idea of AHP is to decompose complex decision problems into multiple levels, establish a hierarchical structure tree, and then allocate quantitative weights to each level, ultimately obtaining the weights of the decision. It can help decision-makers compare and balance multiple decision factors to make optimal decisions\cite{6.5}. However, this method cannot provide new reference solutions for decision-makers, and it is also difficult to determine weights when there are too many factors to evaluate. The ultimate goal of AHP is to fill the weight table through mathematical methods (rather than intuition)\cite{6.6}.

The main steps of Analytic Hierarchy Process (AHP) are as follows\cite{7}:

(1)Establishing a hierarchical structure model: Determine the decision objective and construct a tree-like structure model by organizing the objectives from the overall to the specific level.

(2)Determining the judgment matrix: Compare the importance of each factor pairwise to form a judgment matrix.

(3)Calculating the eigenvectors of the judgment matrix: Calculate the eigenvectors of each judgment matrix and normalize them for weight calculation.

(4)Calculating weights: Calculate the corresponding weights by taking the weighted average of the eigenvectors of each judgment matrix.

(5)Consistency test: Verify the consistency of each judgment matrix to ensure the rationality of weight allocation.

(6)Decision and evaluation: Make decisions and evaluations based on the calculated weights.

\section{Methodology}

In this section, using Analytic Hierarchy Process (AHP), we will establish a judgment matrix between indicators at each level and perform the required consistency check data by calculating the weight index of the indicators\cite{5.8}.

\subsection{Constructing judgment matrix}

When dealing with qualitative issues, it is often necessary to consider multiple factors that are difficult to quantify. However, we still need to know the weight of each factor's impact on our objective. In such cases, the Analytic Hierarchy Process (AHP) can be used to quantify qualitative problems. The ultimate goal of AHP is to mathematically determine the weights, rather than relying on intuition\cite{8}.

If the analytic hierarchy process (AHP) is not used, the weights of actual influencing factors are qualitatively determined, which may be highly unreasonable. Therefore, we need to use a quantitative method (AHP) to characterize this qualitative process, ensuring that it is reasonable and logical. AHP involves comparing the relative importance between pairs of factors and assigning numerical values to indicate their relative importance\cite{8.5}. Typically, experts or decision-makers assess and score the factors based on subjective and objective conditions, using a 1-9 scale. This process quantifies the values of the comparative indicators, with the numerical values reflecting the importance of the indicators, as shown in the Table \ref{Table dingyibiao}.Constructing a Comparative Matrix to Quantify Qualitative Issues.

\begin{table}[!ht]
\centering
\caption{AHP Scale Ranking Definition Table}
\label{Table dingyibiao}
\begin{tabular}{ccc}
\toprule
Scale$a_{ij}$ & Definition  \\
\midrule 
1 & Factor $i$ is equally important as Factor $j$ \\
3 & Factor $i$ is slightly more important than Factor $j$ \\
5 & Factor $i$ is relatively more important than Factor $j$ \\
7 & Factor $i$ is very important compared to Factor $j$ \\
9 & Factor $i$ is of paramount importance compared to Factor $j$ \\
2 4 6 8 & Factor $i$ and Factor $j$ are assigned intermediate-level importance ratings \\
Reciprocal & When comparing Factor $j$ with Factor $i$ and assessing the value as $a_{ji}$,then $a_{ji}$=$\dfrac{1}{a_{ij}}$ \\
\bottomrule
\end{tabular}
\end{table}

The method employed in this study to obtain the judgment matrix involves inviting multiple experts. Through a questionnaire survey, experts provide ratings for the importance of various indicators. The distribution of experts is made as even as possible considering factors such as gender, position, and field. The scoring results are then averaged to obtain the final judgment matrix.

The judgment matrix takes the form:

\begin{equation}
\label{panduanjuzhen}
A=
\begin{pmatrix}
a_{11} & a_{12} & \cdots &  a_{1m} \\
a_{21} & a_{22} & \ldots & a_{2m} \\
\vdots & \vdots & \ddots & \vdots \\
a_{m1} & a_{m2} & \cdots & a_{mm}
\end{pmatrix}
\end{equation}

Each element $a_{ij}$ in matrix A represents the relative importance, with $a_{ji}$=$\dfrac{1}{a_{ij}}$, The elements on the diagonal, $a_{11},a_{22},...a_{mm}$,  are all equal to 1.

\subsection{Calculating weights}

(1)Calculating the Geometric Mean Value of Each Row in the Judgement Matrix(\ref{panduanjuzhen}):

\begin{equation}
\label{pingjunzhi}
\overline{w_{i}}= \sqrt[m] {\prod \limits_{j=1}^ma_{ij}} \quad i=1,2,3...m
\end{equation}

(2)Normalize $\overline{w_{i}}$ :

\begin{equation}
w_{i}= \dfrac{\overline{w_{i}}}{\sum\limits_{i=1}^{m}\overline{w_{i}}} \quad i=1,2,3...m
\end{equation}

Obtain the weight vector $W=\left( w_{1},w_{2},w_{3}...w_{m}\right) ^{T}$, where $w_{i}$ represents the weight of the $i$ indicator.

(3)Calculate the maximum eigenvalue of the judgment matrix \cite{9}:

\begin{equation}
\label{zuidatezhengzhi}
\mu_{\max }= \sum\limits_{i=1}^{m}\dfrac{\left( Aw\right) _{i}}{mw_{i}}
\end{equation}

Where $A$ is the judgment matrix of certain layer indicators obtained, $w$ is the weight vector, and $i$ represents the $i$ th component of the vector.

\subsection{Conducting consistency check}

In the process of determining the weights among various factors at different levels in multi-index comprehensive evaluation, Santy et al. proposed the Judgment Matrix Method. Unlike the approach of comparing all factors together, this method involves pairwise comparisons and utilizes a relative scale to mitigate the difficulty of comparing factors with different natures, thereby enhancing the accuracy of the comparisons.

However, during the process of pairwise comparisons, achieving complete consistency in judgments is not feasible, giving rise to the issue of estimation errors. This inevitably leads to biases in eigenvalues and weight vectors, and can also result in contradictory situations, which are objectively present. To prevent significant errors, it is necessary to assess the consistency of the judgment matrix\cite{10}.

Consistency check is performed to validate the harmony among the importance of multiple elements, avoiding contradictions such as A being more important than B, B being more important than C, and C being more important than A. When conducting a comprehensive multi-index evaluation, using a consistency check can effectively assess the degree of consistency in the judgment matrix, thus ensuring the accuracy of weight calculations.

The steps for conducting a consistency check are as follows:

(1)Construct the consistency indicators:

\begin{equation}
\label{CI}
CI=\dfrac{\mu _{\max }-m}{m-1}
\end{equation}

Where $m$ is the order of the matrix.When $CI=0$, the matrix achieves complete consistency.The larger the $CI$, the lower the matrix consistency. From the formula (\ref{CI}) , it is evident that the value of $CI$ is related to the order $m$ of the matrix. To mitigate the adverse impact of the order $m$, an adjustment should be made by introducing a random consistency index $RI$.

(2)Refer to the literature to introduce $RI$, the Random Index for Consistency, as shown in the Table \ref{Table RI}.This value serves as a corrective coefficient for matrix inconsistency and acts as a remedy for the influence of matrix order. Table \ref{Table RI} represents Satty's calculations for different $m$, derived from 1000 sample matrices $A_{1}$: For a fixed $m$, random positive reciprocal matrices $A_{1}$ are generated, and then the consistency index $CI=\dfrac{\mu _{\max }-m}{m-1}$ is computed for each $A_{1}$. $A_{1}$ is highly inconsistent, resulting in significantly large $CI$ values. By constructing a substantial number of $A_{1}$ matrices and averaging their $CI$ values, the random consistency index $RI$ is determined.

\begin{table}[!ht]
\centering
\caption{Random Index $RI$ Value}
\label{Table RI}
\begin{tabular}{cccccccccccc}
\toprule
m & 1 & 2 & 3 & 4 & 5 & 6 & 7 & 8 & 9 & 10 & 11 \\
\midrule 
$RI$ & 0 & 0 & 0.58 & 0.90 & 1.12 & 1.24 & 1.32 & 1.41 & 1.45 & 1.49 & 1.51 \\
\bottomrule 
\end{tabular}
\end{table}

(3)Construct Consistency Ratio:

\begin{equation}
CR=\dfrac{CI}{RI}
\end{equation}

When $CR<0.1$, it can be considered that the inconsistency level of the matrix is within an acceptable range, and the matrix passes the consistency check.

\section{Application}

In this section, using the relevant concepts of AHP and specific obtained judgment matrices, we can not only obtain the weights of all energy efficiency indicators required by some cigarette manufacturing enterprise, but also analyze their weight ranking.

\subsection{Obtain the judgment matrix}

The establishment of the judgment matrix and its distribution of weights in Analytic Hierarchy Process (AHP) is based on the pairwise comparisons of the relative importance of factors at each level with respect to the corresponding factors at the higher level. This task was accomplished by sending survey questionnaires to experts within a certain cigarette enterprise. Quantitative judgment matrices such as $A-B$ and $B-C$ were formed by multiple experts assigning scores using a 1-9 proportional scale, and then the average values were calculated to derive the judgment matrices used for computation. The final outcome is as follows: (\ref{A-B})-(\ref{B_{6}-C}).

\begin{equation}
\label{A-B}
A=
\begin{pmatrix}
1 & 1.3803 & 1.5556 & 1.5806 &  1.5077 & 2.3902 \\
0.7245 & 1 & 1.1270 & 1.1452 &  1.0923 & 1.7317 \\
0.6429 & 0.8873 & 1 & 1.0161 &  0.9692 & 1.7317 \\
0.6327 & 0.8732 & 0.9841 & 1 &  0.9538 & 1.5122 \\
0.6633 & 0.9155 & 1.0317 & 1.0484 &  1 & 1.5854 \\
0.4184 & 0.5775 & 0.6508 & 0.6613 &  0.6308 & 1 \\
\end{pmatrix}
\end{equation}

\begin{equation}
\label{B_{2}-C}
B_{2}-C=
\begin{pmatrix}
1 & 1.4265 & 0.6510 & 0.5215  \\
0.7010 & 1 & 0.4564 & 0.3656  \\
1.5361 & 2.1912 & 1 & 0.8011  \\
1.9175 & 2.7353 & 1.2483 & 1  \\
\end{pmatrix}
\end{equation}

\begin{equation}
\label{B_{3}-C}
B_{3}-C=
\begin{pmatrix}
1 & 1.5952 & 0.9054 & 0.7882 &  1.0806 & 0.9571 \\
0.6269 & 1 & 0.5676 & 0.4941 &  0.6774 & 0.6000 \\
1.1045 & 1.7619 & 1 & 0.8706 &  1.1935 & 1.0571 \\
1.2687 & 2.0238 & 1.1486 & 1 &  1.3710 & 1.2143 \\
0.9254 & 1.4762 & 0.8378 & 0.7294 &  1 & 0.8857 \\
1.0448 & 1.6667 & 0.9459 & 0.8235 &  1.1290 & 1 \\
\end{pmatrix}
\end{equation}

\begin{equation}
\label{B_{4}-C}
B_{4}-C=
\begin{pmatrix}
1 & 1.6100 & 2.0214 & 1.9389  \\
0.6211 & 1 & 1.2553 & 1.2041  \\
0.4947 & 0.7966 & 1 & 0.9592  \\
0.5158 & 0.8305 & 1.0426 & 1  \\
\end{pmatrix}
\end{equation}

\begin{equation}
\label{B_{5}-C}
B_{5}-C=
\begin{pmatrix}
1 & 1.8387 & 1.6765 & 1.6474 &  1.3103  \\
0.5439 & 1 & 0.9118 & 0.8986 &  0.7126  \\
0.5965 & 1.0968 & 1 & 0.9855 &  0.7816  \\
0.6053 & 1.1129 & 1.0147 & 1 &  0.7931 \\
0.7632 & 1.4032 & 1.2794 & 1.2609 &  1  \\
\end{pmatrix}
\end{equation}

\begin{equation}
\label{B_{6}-C}
B_{6}-C=
\begin{pmatrix}
1 & 1.1153 & 1.7652    \\
0.8966 & 1 & 1.5827    \\
0.5665 & 0.6318 & 1    \\
\end{pmatrix}
\end{equation}

\subsection{Calculating the weights of various indicators in the cigarette manufacturing enterprise's production efficiency system}

Calculate each indicator in layer $C$ according to Formula (\ref{pingjunzhi})-(\ref{zuidatezhengzhi}), and compare it with the weights of layer $B$; calculate each indicator in layer $B$ and compare it with the weights of layer $A$. Verify the consistency between the two sets of comparisons, and summarize the results in Table \ref{Table yizhixingjianyan}.

\begin{table}[!ht]
\centering
\caption{Indicator Weights and Consistency Test Results}
\label{Table yizhixingjianyan}
\resizebox{\linewidth}{!}{
\begin{tabular}{cccccccc}
\toprule
Matrix & Normalized Weight Vectors $W$ & Maximum Eigenvalue $\mu _{\max }$ & $CI$ & $RI$ & $CR$ & Consistency Test Results \\
\midrule 
$A-B$   & $\omega= [0.245 , 0.1775 , 0.1575, 0.155, 0.1625, 0.1025]$   &  6   & 0 & 1.25	& 0   &  Passed     \\
$B_{1}-C$   &  $\omega= [1]$   &  1  & 0 & 0	& 0  &    Passed   \\
$B_{2}-C$   &  $\omega= [0.194 , 0.136 , 0.298 , 0.372]$   & 4	& 0	& 0.882	& 0   &   Passed    \\
$B_{3}-C$   &  $\omega= [0.1675 , 0.105 , 0.185 , 0.2125 , 0.155 , 0.175]$   & 6 & 0 & 1.25	& 0   &    Passed   \\
$B_{4}-C$   &  $\omega= [0.38 , 0.236 , 0.188 , 0.196]$    &  4	& 0	& 0.882	& 0   &    Passed   \\
$B_{5}-C$   &  $\omega= [0.285 , 0.155 , 0.17 , 0.1725 , 0.2175]$    & 5	& 0	& 1.11	& 0        &    Passed   \\
$B_{6}-C$   &  $\omega= [0.406 , 0.364 , 0.23]$   &  3	& 0	& 0.525	& 0   &    Passed   \\
\bottomrule 
\end{tabular}
}
\end{table}

\subsection{Weight ranking}

The Table \ref{Table quanzhonghuizong} below summarizes the weightings of performance evaluation indicators for some Cigarette Factory, based on the comparison of weights between indicator layers $C$ and $B$.

\begin{table}[!ht]
\centering
\caption{Ranking of Comprehensive Weights for A Certain Cigarette Factory Indicators}
\label{Table quanzhonghuizong}
\resizebox{\linewidth}{!}{
\begin{tabular}{ccccccccccc}
\toprule
\multirow{2}{*}{Weight} & $B_{1}$ & $B_{2}$ & $B_{3}$ & $B_{4}$ & $B_{5}$ & $B_{6}$  & \multirow{2}{*}{Overall Weight} \\
                    & 0.245  & 0.1775  & 0.1575  & 0.155  &  0.1625 & 0.1025    &                    \\ \hline
$C_{11}$Production Plan Completion Rate                   & 0.245  &   &   &   &   &   &      0.245                    \\
$C_{24}$Packaging Equipment Operating Efficiency                   &   &  0.372 &   &   &   &   &      0.06603                    \\
$C_{41}$Number of Process Deviations                   &   &   &   & 0.38  &   &   &      0.0589                    \\
$C_{23}$Rolling Equipment Operating Efficiency                   &   &  0.298 &   &   &   &   &      0.052895                    \\
$C_{51}$$CPK$Qualification Rate of Cigarette Quality                   &   &   &   &   & 0.285  &   &      0.0463125                    \\
$C_{61}$Energy Consumption Per Carton                   &   &   &   &   &   & 0.406  &      0.041615                    \\
$C_{62}$Pollution Emissions Per Carton                   &   &   &   &   &   & 0.364  &      0.03731                    \\
$C_{42}$Acceptance Rate of Premium Quality Finished Products                   &   &   &   & 0.236  &   &   &      0.03658                    \\
$C_{55}$Premium Product Rate                   &   &   &   &   & 0.2175  &   &      0.03534375                    \\
$C_{21}$Unit Output Maintenance Duration                   &   &  0.194 &   &   &   &   &      0.034435                    \\
$C_{34}$Residual Cigarettes Per Carton                   &   &   & 0.2125  &   &   &   &      0.03346875                    \\
$C_{44}$Moisture Deviation in Fuel Drying Machine Export                   &   &   &   & 0.196  &   &   &      0.03038                    \\
$C_{43}$Moisture Deviation in Thin Sheet Drying Machine Export                   &   &   &   & 0.188  &   &   &      0.02914                    \\
$C_{33}$Consumption of Small Box Paper Per Carton                  &   &   & 0.185  &   &   &   &       0.0291375                   \\
$C_{54}$Score of Finished Product Release Inspection                   &   &   &   &   & 0.1725  &   &      0.02803125                    \\
$C_{53}$Deviation in Single Cigarette Weight                   &   &   &   &   & 0.17  &   &      0.027625                    \\
$C_{36}$Consumption of Cut Cigarette Per Carton                   &   &   & 0.175  &   &   &   &      0.0275625                    \\
$C_{31}$Leaf-to-Silk Ratio                   &   &   & 0.1675  &   &   &   &      0.02638125                    \\
$C_{52}$$CPK$Qualification Rate of Draw Resistance Quality                  &   &   &   &   & 0.155  &   &      0.0251875                    \\
$C_{35}$Consumption of Cigarette Filters Per Carton                   &   &   & 0.155  &   &   &   &      0.0244125                    \\
$C_{22}$Unit Output Maintenance Frequency                   &   &  0.136 &   &   &   &   &      0.02414                    \\
$C_{63}$Carbon Dioxide Emissions Per Carton                   &   &   &   &   &   & 0.23  &      0.023575                    \\
$C_{32}$Consumption of Cigarette Paper Per Carton                   &   &   & 0.105  &   &   &   &      0.0165375                    \\
\bottomrule 
\end{tabular}
}
\end{table}

\section{Model Analysis of Cigarette Manufacturing Company's Production Efficiency Based on Analytic Hierarchy Process}

This section will conduct an analysis of the production efficiency model of a certain cigarette production enterprises based on the AHP, mainly selecting the required indicators based on importance ranking.

\subsection{Comparison of Weights between Layer $B$ Indicators and Layer $A$ Indicators}

The calculated weights of Layer $B$ indicators relative to Layer $A$ are presented in Table \ref{Table quanzhonghuizong}. Based on the numerical values of the weights, the factors within Criterion Layer $B$ are ranked as follows:

$\omega _{B_{1}}=0.245 > \omega _{B_{2}}=0.1775 > \omega _{B_{5}}=0.1625 > \omega _{B_{3}}=0.1575 > \omega _{B_{4}}=0.155 > \omega _{B_{6}}=0.1025$

Among these indicators, the weight of ${B_{1}}$ is greater than 0.2, indicating that 'Plan Execution' plays a significant role in Criterion Layer $B$ and is the most important indicator within Layer $B$. The weights of ${B_{2}}$ ${B_{5}}$ ${B_{3}}$ ${B_{4}}$ are relatively close, with a difference of only 0.02, suggesting that 'Equipment Efficiency' 'Process Control of Rolling and Packaging Quality' 'Production Material Consumption' and 'Process Control of Silk Production Quality' are nearly equally important. The weight of ${B_{6}}$ is the lowest, but this doesn't imply that "Energy Saving and Emission Reduction" is unimportant; it simply means that, in the efficiency system of a certain cigarette enterprise, it is relatively less prioritized compared to the first five indicators.

\subsection{Comparison of Weights between Layer $C$ Indicators and Layer $B$ Indicators}

The calculated weights of Layer $C$ indicators relative to Layer $B$ are presented in Table \ref{Table quanzhonghuizong}. Based on the numerical values of the weights, the factors within Criterion Layer $C$ are ranked as follows:

(1)The weight rankings for the three-level indicators under Criterion ${B_{2}}$ are as follows:

$\omega _{C_{24}}=0.372 > \omega _{C_{23}}=0.298 > \omega _{C_{21}}=0.194 > \omega _{C_{22}}=0.136$

Among these indicators, both ${C_{24}}$ and ${C_{23}}$ have weights greater than 0.29, and their sum is 0.67, close to two-thirds of the total weights. Therefore, ``Packaging Equipment Operating Efficiency'' and ``Rolling Equipment Operating Efficiency'' play a significant role in criterion layer ${B_{2}}$ and are the two most important indicators in ${B_{2}}$. Although ${C_{21}}$ and ${C_{22}}$ have relatively lower rankings in terms of weight, it does not mean that the indicators ``Unit Output Maintenance Duration'' and ``Unit Output Maintenance Frequency'' are unimportant. It simply indicates that they have a relatively lower impact compared to the previous two indicators in the "Equipment Efficiency" criterion layer of ${B_{2}}$.

(2)The weight rankings for the three-level indicators under Criterion ${B_{3}}$ are as follows:

$\omega _{C_{34}}=0.2125 > \omega _{C_{33}}=0.185 > \omega _{C_{36}}=0.175 > \omega _{C_{31}}=0.1675 > \omega _{C_{35}}=0.155 > \omega _{C_{32}}=0.105$

Among these indicators, the weight of indicator ${C_{34}}$ is greater than 0.2, so "Residual Cigarettes Per Carton" plays a significant role in criterion layer ${B_{3}}$ and is the most important indicator in ${B_{3}}$.The weights of the four indicators ${C_ {33}}$, ${C_{36}}$, ${C_{31}}$ and ${C_{35}}$ are relatively close, with a difference of only 0.03, indicating that "Consumption of Small Box Paper Per Carton", "Consumption of Cut Cigarette Per Carton", "Leaf-to-Silk Ratio", and "Consumption of Cigarette Filters Per Carton" are nearly equally important. The weight of indicator ${C_{32}}$ is the lowest, but it doesn't mean that the indicator "Consumption of Cigarette Paper Per Carton" is unimportant. It is just relatively less considered compared to the previous five indicators in criterion layer ${B_{3}}$ "Production Material Consumption".

(3)The weight rankings for the three-level indicators under Criterion ${B_{4}}$ are as follows:

$\omega _{C_{41}}=0.38 > \omega _{C_{42}}=0.236 > \omega _{C_{44}}=0.196 > \omega _{C_{43}}=0.188$

The weights of indicators ${C_{41}}$ and ${C_{42}}$ are both greater than 0.23, and their sum is 0.61, which is close to two-thirds of the total weight. Therefore, "Number of Process Deviations" and "Acceptance Rate of Premium Quality Finished Products" play a significant role in criterion layer ${B_{4}}$ and are the two most important indicators in ${B_{4}}$. Although indicators ${C_{44}}$ and ${C_{43}}$ have lower weight rankings, it does not mean that "Moisture Deviation in Fuel Drying Machine Export" and "Moisture Deviation in Thin Sheet Drying Machine Export" are unimportant. They are relatively less influential compared to the previous two indicators in the "Process Control of Silk Production Quality" criterion layer in ${B_{4}}$.

(4)The weight rankings for the three-level indicators under Criterion ${B_{5}}$ are as follows:

$\omega _{C_{51}}=0.285 > \omega _{C_{55}}=0.2175 > \omega _{C_{54}}=0.1725 > \omega _{C_{53}}=0.17 > \omega _{C_{52}}=0.155$

The weights of both indicators ${C_{51}}$ and ${C_{55}}$ are greater than 0.2, and their sum is 0.5, which accounts for half of the total weight. Therefore, "$CPK$Qualification Rate of Cigarette Quality" and "Premium Product Rate" play a significant role in criterion layer ${B_{5}}$ and are the two most important indicators in ${B_{5}}$. Although the weights of indicators ${C_{54}}$, ${C_{53}}$ and ${C_{52}}$ are ranked lower, it does not mean that "Score of Finished Product Release Inspection", "Deviation in Single Cigarette Weight" and "$CPK$Qualification Rate of Draw Resistance Quality" are not important. These indicators are relatively less influential in the "Process Control of Rolling and Packaging Quality" criterion layer in ${B_{5}}$.

(5)The weight rankings for the three-level indicators under Criterion ${B_{6}}$ are as follows:

$\omega _{C_{61}}=0.406 > \omega _{C_{62}}=0.364 > \omega _{C_{63}}=0.23$

The weights of both indicators ${C_{61}}$ and ${C_{62}}$ are greater than 0.3, and their sum is 0.78. Therefore, "Energy Consumption Per Carton" and "Pollution Emissions Per Carton" play a significant role in criterion layer ${B_{6}}$ and are the two most important indicators in ${B_{6}}$. Although the weight of indicators ${C_{63}}$ is lower, it does not mean that "Carbon Dioxide Emissions Per Carton" is not important. It is relatively less influential compared to the previous two indicators in the "Carbon Dioxide Emissions Per Carton" criterion layer in ${B_{6}}$.

\subsection{Comparison of Weights between Layer $C$ Indicators and Layer $A$ Indicators}

By calculation, the comparison results of the weights between the indicators in layer $C$ and layer $A$ are obtained. Sorting the weight results, the summarized sorting results are shown in Table \ref{Table quanzhonghuizong}. It can be observed that the sorting results are as follows:

$\omega _{C_{11}}=0.245 >\omega _{C_{24}}=0.06603 > \omega _{C_{41}}=0.0589 >\omega _{C_{23}}=0.052895 > \omega _{C_{51}}=0.0463125 > \omega _{C_{61}}=0.041615 > \omega _{C_{62}}=0.03731 > \omega _{C_{42}}=0.03658 > \omega _{C_{55}}=0.03534375 > \omega _{C_{21}}=0.034435 > \omega _{C_{34}}=0.03346875 > \omega _{C_{44}}=0.03038 > \omega _{C_{43}}=0.02914 > \omega _{C_{33}}=0.0291375 > \omega _{C_{54}}=0.02803125 > \omega _{C_{53}}=0.027625 > \omega _{C_{36}}=0.0275625 > \omega _{C_{31}}=0.02638125 > \omega _{C_{52}}=0.0251875 > \omega _{C_{35}}=0.0244125 > \omega _{C_{22}}=0.02414 > \omega _{C_{63}}=0.023575 > \omega _{C_{32}}=0.0165375$

The weight of $C_{11}$ is the highest, being 0.245, which is close to one-fourth of the total weight. This indicates that among the 23 three-level indicators, "Production Plan Completion Rate" is the most important factor. Whether the monthly production plan is completed directly affects the efficiency evaluation of the cigarette factory for that month.

The weights of ${C_{24}}$, ${C_{41}}$, ${C_{23}}$, ${C_{51}}$ and $C_{61}$ are all greater than 0.04. They are among the most important indicators in the target layer $A$. Their total weight is 0.26, which also reaches one-fourth of the total weight. Adding the weight of $C_{61}$, their combined weight accounts for half of the total weight. This indicates that "Production Plan Completion Rate", "Packaging Equipment Operating Efficiency", "Number of Process Deviations", "Rolling Equipment Operating Efficiency", "$CPK$Qualification Rate of Cigarette Quality", and "Energy Consumption Per Carton" are the most important 6 indicators compared to the remaining 17 tertiary indicators.

The weights of ${C_{62}}$ "Pollution Emissions Per Carton", ${C_{42}}$ "Acceptance Rate of Premium Quality Finished Products", ${C_{55}}$ "Premium Product Rate", ${C_{21}}$ "Unit Output Maintenance Duration", ${C_{34}}$ "Residual Cigarettes Per Carton" and ${C_{44}}$ "Moisture Deviation in Fuel Drying Machine Export" are all greater than 0.03, their total weight is 0.2.

The weights of ${C_{43}}$ "Moisture Deviation in Thin Sheet Drying Machine Export", ${C_{33}}$ "Consumption of Small Box Paper Per Carton", ${C_{54}}$ "Score of Finished Product Release Inspection", ${C_{53}}$ "Deviation in Single Cigarette Weight", ${C_{36}}$ "Consumption of Cut Cigarette Per Carton", ${C_{31}}$ "Leaf-to-Silk Ratio" and ${C_{52}}$ "$CPK$Qualification Rate of Draw Resistance Quality" are all greater than 0.025, their total weight is 0.19.

The four indicators ${C_{35}}$ "Consumption of Cigarette Filters Per Carton", ${C_{22}}$ "Unit Output Maintenance Frequency", ${C_{63}}$ "Carbon Dioxide Emissions Per Carton" and ${C_{32}}$ "Consumption of Cigarette Paper Per Carton" have weights ranging from 0.016 to 0.025. The total weight of these indicators is 0.09, indicating their relatively lower importance.

\section{Conclusion}

This paper has constructed a more scientific, objective, comprehensive, and universal performance indicator system for cigarette enterprises. The complex evaluation system is hierarchically structured, and quantitative weight allocation is applied to each level to obtain the final decision weights. This indicator system can assist cigarette enterprises in comparing and balancing between various indicators, facilitating a more intuitive analysis of problems.

Based on the constructed performance indicator system for cigarette enterprises, the opinions of cigarette experts were obtained through the Delphi method and applied in practical scenarios to derive the final indicator weights. Through analysis, it becomes evident which production indicators are most important in assessing production efficiency, providing reference for future operations of cigarette enterprises.

\section{Data Availability}

\hspace{2em}
The data in this article has been provided or preliminarily summarized and analyzed by the first, third, fourth, and fifth authors. The datasets generated during and/or analysed during the current study are available from the corresponding author on reasonable request.

\section{Declarations}
\hspace{2em}The authors declare that there are no conflict of interests, we do not have any possible conflicts of interest.

\vskip 6mm
\noindent{\bf Acknowledgements}

\noindent
\hspace{2em}The third author is partly supported by the China Tobacco Hebei industrial co.,Ltd Technology Project  under Grant No. HBZY2023A034, and the sixth author is partly supported by the Tianjin Natural Science Foundation under Grant No.18JCYBJC16300.


\begin{thebibliography}{99}
\bibitem{0.5}M.Li , Development of Quality Control Index System and Supporting Detection Techniques for Water-Based Adhesives Used in Tobacco Industry, D. Hunan University, 2019.
\bibitem{0.6}F.Hu , Digital Evaluation System for the Tobacco Leaf Addition Process and Research on Addition Effects Prediction, D. Kunming University of Science and Technology, 2020.
\bibitem{1}S.F.Li , Research on the Construction of Lean Production Management System in the Cigarette Industry Enterprises, Technological Innovation and Productivity. 03 (2023), 78-80.
\bibitem{1.4}X.R.Chen , Research on Evaluation System for Suppliers of Tobacco Additives, D. Kunming University of Science and Technology, 2007.
\bibitem{1.5}G.Y.Liao , Construction of Smoking Control Effect Evaluation Index System in Public Places in Chongqing City, J. Modern Preventive Medicine, 2023.
\bibitem{2}X.Y.Chen , Exploring the Cost Control Model of Lean Production in Cigarette Factories, Money China. 01 (2022), 55-57.
\bibitem{3}B.Wang , Comprehensive evaluation of urban garden afforestation based on PLS-SEM path, J. Physics and Chemistry of the Earth,2022, 126.
\bibitem{3.5}Y.W.Gao , Construction of Measurable Impact Evaluation Indicator System for Open Science Datasets, J. Information Science, 07(2023), 49-60.
\bibitem{4}Qian Z , Getu N ,Wei J , et al, Research on key index evaluation of power transmission and transformation wiring based on three-dimensional intelligent evaluation,J. Energy Reports,2022, 8(S7).
\bibitem{4.5}X.R.Sun ,J.H.Sun , Construction of Elderly Care Assessment Indicator System Based on Basic Abilities and Basic States Using Delphi Method and Analytic Hierarchy Process, J. Journal of Practical Cardiocerebral Vascular Disease, 02(2018), 58-62.
\bibitem{5}J.M.Niu ,D.Wang ,X.Li , et al, Construction of a fuzzy comprehensive evaluation model for the evaluation of cigarette specification competitiveness in tobacco commercial enterprises, J. Market Modernization, 06 (2022), 1-5.
\bibitem{5.8}P.P.Zhang ,Z.R.Li ,H.Y.Song , et al, Research on Site Selection Indicator System for Deep-Water Cage Aquaculture Based on AHP, J. Southern Fisheries Science, 04(2023), 1-9.
\bibitem{6}Saaty, Thomas L. Deriving the AHP 1-9 scale from first principles. Proceedings 6th ISAHP. Berna, Suiza (2001), 397-402.
\bibitem{6.5}B.Xv , Research on the Construction of Undergraduate Classroom Teaching Quality Evaluation Indicator System with a Student-Centered Approach, D. Northeast Petroleum University, 2023.
\bibitem{6.6}J.X.Song ,K.X.Wei , Research on River Sustainability Assessment Indicator System and Methods, J. Yangtze River, 05(2023), 40-46.
\bibitem{7}Luo , Tingyu, et al , Comprehensive Evaluation of Emergency Shelters in Wuhan City Based on GIS. 2022 29th International Conference on Geoinformatics. IEEE, (2022), 1-6.
\bibitem{8}J.J.Liu , Research on the Value of the ``Li Canal - Gaoyou Irrigation District'' Based on AHP-Fuzzy Comprehensive Evaluation Method,D. Anhui University of Architecture, 2023.
\bibitem{8.5}Y.Peng ,T.Xv , Research on the Construction of Agricultural Product Logistics Evaluation System Based on AHP, J. Chinese Business Review, 14(2023), 101-104.
\bibitem{9}Y.S.Yang , Teaching Difficulties in Calculating Indicator Weights Using Analytic Hierarchy Process, J. Science Education and Culture, 05(2018), 67-68.
\bibitem{10}Y.H.Tong ,Z.L.Feng , et al, Evaluation and Analysis of Haze Influencing Factors Based on AHP, J.  Journal of Southwest Normal University (Natural Science Edition), 03(2020), 87-94.



\end{thebibliography}
\end{document}